\RequirePackage{ifpdf}
\ifpdf % We are running pdfTeX in pdf mode
\documentclass[pdftex]{sigma}
\else
\documentclass{sigma}
\fi

\begin{document}

\renewcommand{\PaperNumber}{113}

\FirstPageHeading

\renewcommand{\thefootnote}{$\star$}

\ShortArticleName{Symmetries and Invariant Dif\/ferential Pairings}

\ArticleName{Symmetries and Invariant Dif\/ferential Pairings\footnote{This paper is a contribution to the Proceedings
of the Seventh International Conference ``Symmetry in Nonlinear
Mathematical Physics'' (June 24--30, 2007, Kyiv, Ukraine). The
full collection is available at
\href{http://www.emis.de/journals/SIGMA/symmetry2007.html}{http://www.emis.de/journals/SIGMA/symmetry2007.html}}}

\Author{Michael G. EASTWOOD}

\AuthorNameForHeading{M.G. Eastwood}

\Address{Department of Mathematics, University of Adelaide,
SA 5005, Australia}

\Email{\href{mailto:meastwoo@member.ams.org}{meastwoo@member.ams.org}}

\ArticleDates{Received November 14, 2007; Published online November 23, 2007}

\Abstract{The purpose of this article is to motivate the study of invariant,
and especially conformally invariant, dif\/ferential pairings. Since a general
theory is lacking, this work merely presents some interesting examples of
these pairings, explains how they naturally arise, and formulates various
associated problems.}


\Classification{53A30; 58J70; 53A20}

\section{Introduction}\label{intro}
The Lie derivative is an extremely familiar operation in dif\/ferential geometry.
Given a smooth vector f\/ield $V$ and tensor f\/ield~$\phi$, the Lie derivative
${\mathcal{L}}_V\phi$ is certainly an intrinsic or invariant construction. In
practise, we can take this to mean that the result is independent of any local
co\"ordinate formula that may be used. Alternatively, we may employ an
arbitrary torsion-free af\/f\/ine connection $\nabla$ to write explicit formulae
such as
\begin{gather}\label{Lie}
\phi_{bc}\stackrel{{\mathcal{L}}_V}{\longmapsto}
V^a\nabla_a\phi_{bc}+(\nabla_bV^a)\phi_{ac}+(\nabla_cV^a)\phi_{ba},
\end{gather}
checking that the result is independent of the connection used. In this
formula, the indices are `abstract indices' in the sense of Penrose~\cite{OT}.
It follows that the expression is co\"ordinate-free whilst its freedom from
choice of connection follows immediately from the formulae
\begin{gather}\label{contorsion}
\widehat\nabla_aV^c=\nabla_aV^c+\Gamma_{ab}{}^cV^b\qquad\mbox{and}\qquad
\widehat\nabla_a\phi_{bc}=\nabla_a\phi_{bc}-\Gamma_{ab}{}^d\phi_{dc}
-\Gamma_{ac}{}^d\phi_{bd},\end{gather}
capturing all possible torsion-free connections via the choice of tensor
$\Gamma_{ab}{}^c=\Gamma_{ba}{}^c$. For further details on this point of view,
see~\cite{OT}.

Writing the Lie derivative as ${\mathcal{L}}_V\phi$ suggests that we are
regarding it as a linear dif\/ferential operator acting on the tensor f\/ield
$\phi$ for a f\/ixed choice of vector f\/ield~$V$. However, it is clear from
expressions such as (\ref{Lie}) that we may equally well f\/ix $\phi$ and regard
the result as a linear dif\/ferential operator acting on vector f\/ields~$V$.
Sometimes, this is the natural viewpoint. For example, if $\phi_{ab}$ is
symmetric then we may rewrite (\ref{Lie}) as
\[({\mathcal{L}}_V\phi)_{bc}=\phi_{bc}\nabla_{a}V^c+\phi_{ac}\nabla_{b}V^c
+(\nabla_c\phi_{ab})V^c\]
and now if $\phi_{ab}$ is a metric and we choose $\nabla_a$
to be its Levi-Civita connection, then
\[{\mathcal{L}}_V\phi=0 \ \iff \ \nabla_aV_b+\nabla_bV_a=0,\]
which is the usual way of viewing the Killing equation on a Riemannian
manifold. Of course, the best way of viewing ${\mathcal{L}}_V\phi$ is as a
{\em bilinear} dif\/ferential operator or as a {\em differential pairing}. It is
invariantly def\/ined on any manifold and is an example of the subject matter of
this article.

One can easily imagine other expressions along the lines of the right hand side
of (\ref{Lie}) that turn out to be independent of choice of connection.
Restricting the choice of connection leads to a greater range of expressions
that might be invariant. In the next two sections we shall f\/ind examples that
are invariant when the connections are deemed to be Levi-Civita connections for
metrics taken from a f\/ixed conformal class.

\section{First order symmetries}\label{symmetries}
Suppose $\mbox{\L}:E\to F$ is a linear dif\/ferential operator between vector
bundles $E$ and~$F$. A linear dif\/ferential operator ${\mathcal{D}}:E\to E$ is
said to be a {\em symmetry} of $\mbox{\L}$ if and only if
\[\mbox{\L}{\mathcal{D}}=\delta \mbox{\L}\quad
\mbox{for some linear dif\/ferential operator }\delta:F\to F.\]
When $\mbox{\L}$ is the Laplacian on ${\mathbb{R}}^3$ (and the bundles $E$ and
$F$ are trivial), its f\/irst order symmetries are well-known~\cite{miller}.
Specif\/ically, they comprise an $11$-dimensional vector space spanned by the
following
\begin{gather}
f\mapsto f,\qquad
f\mapsto\frac{\partial f}{\partial x},\qquad
f\mapsto\frac{\partial f}{\partial y},\qquad
f\mapsto\frac{\partial f}{\partial z},\qquad
f\mapsto x\frac{\partial f}{\partial x}+y\frac{\partial f}{\partial y}
+z\frac{\partial f}{\partial z},\nonumber\\
f\mapsto x\frac{\partial f}{\partial y}-y\frac{\partial f}{\partial x},\qquad
f\mapsto (x^2-y^2-z^2)\frac{\partial f}{\partial x}
+2xy\frac{\partial f}{\partial y}+2xz\frac{\partial f}{\partial z}+xf,\nonumber\\
f\mapsto y\frac{\partial f}{\partial z}-z\frac{\partial f}{\partial y},\qquad
f\mapsto (y^2-z^2-x^2)\frac{\partial f}{\partial y}
+2yz\frac{\partial f}{\partial z}+2yx\frac{\partial f}{\partial x}+yf,\nonumber\\
f\mapsto z\frac{\partial f}{\partial x}-x\frac{\partial f}{\partial z},\qquad
f\mapsto (z^2-x^2-y^2)\frac{\partial f}{\partial z}
+2zx\frac{\partial f}{\partial x}+2zy\frac{\partial f}{\partial y}+zf,\label{eleven}
\end{gather}
where $x$, $y$, $z$ are the usual Euclidean co\"ordinates on~${\mathbb{R}}^3$.
Following~\cite{laplace}, it is convenient to rewrite the general f\/irst order
symmetry as
\[  f\mapsto V^a\nabla_af+\tfrac{1}{6}(\nabla_aV^a)f+Cf,\]
where $V^a$ is an arbitrary vector f\/ield of the form
\begin{gather}\label{confK}
 V^a=-s^a-m^a{}_bx^b+\lambda x^a+r_bx^bx^a-\tfrac{1}{2}x_bx^br^a
\end{gather}
and $C$ is an arbitrary constant. In fact, vector f\/ields of the form
(\ref{confK}) on ${\mathbb{R}}^n$ are precisely the conformal Killing f\/ields,
i.e.~the solutions of the equation
\begin{gather}\label{cKillingfield}
\textstyle\nabla_aV_b+\nabla_bV_a=\frac{2}{n}g_{ab}\nabla^cV_c,\end{gather}
where $g_{ab}$ is the usual Euclidean metric on~${\mathbb{R}}^n$. More
generally, on ${\mathbb{R}}^n$ the f\/irst order symmetries of the Laplacian may
be written as $f\mapsto{\mathcal{D}}_Vf+Cf$, where
\begin{gather}\label{DVf}
{\mathcal{D}}_Vf\equiv V^a\nabla_af+\tfrac{n-2}{2n}(\nabla_aV^a)f
\end{gather}
for an arbitrary conformal Killing f\/ield~$V^a$. This expression certainly
resembles the Lie derivative and we now make this precise as follows. If
$\phi_{bc\cdots de}$ is an $n$-form, then
\[
0=(n+1)(\nabla_{[a}V^a)\phi_{bc\cdots de]}=
(\nabla_aV^a)\phi_{bc\cdots de}-\big[(\nabla_bV^a)\phi_{ac\cdots de}
+\cdots+(\nabla_eV^a)\phi_{bc\cdots da}\big],
\]
where we are following \cite{OT} in employing square brackets to denote the
result of skewing over the indices they enclose and also noting that a skew
tensor with $n+1$ indices on ${\mathbb{R}}^n$ necessarily vanishes.
Hence, the Lie derivative on $n$-forms simplif\/ies:
\begin{gather*}
{\mathcal{L}}_V\phi_{bc\cdots de}=
V^a\nabla_a\phi_{bc\cdots de}+(\nabla_bV^a)\phi_{ac\cdots de}
+\cdots+(\nabla_eV^a)\phi_{bc\cdots da}\\
\phantom{{\mathcal{L}}_V\phi_{bc\cdots de}}{} =V^a\nabla_a\phi_{bc\cdots de}+(\nabla_aV^a)\phi_{bc\cdots de}.
\end{gather*}
Thus, if $h$ is a section
of the bundle $\Lambda^n$ of $n$-forms, then
\[
{\mathcal{L}}_Vh=V^a\nabla_ah+(\nabla_aV^a)h.
\]
Therefore, the expression (\ref{DVf}) coincides with the Lie derivative if $f$
is interpreted as a section of the line bundle $(\Lambda^n)^{(n-2)/2n}$
(supposing that the manifold is orientable so that such fractional powers are
allowed). In particular, we deduce that
$[{\mathcal{D}}_V,{\mathcal{D}}_W]={\mathcal{D}}_{[V,W]}$, as is also readily
verif\/ied by direct computation. Also we conclude that, with this
interpretation, the dif\/ferential pairing~${\mathcal{D}}_Vf$ between $V^a$ and
$f$ is invariant in the sense discussed for Lie derivative in~Section~\ref{intro}.

\section{Conformal geometry}\label{conformalgeometry}

An alternative viewpoint on these matters is as follows. If $\nabla_a$ is the
Levi-Civita connection for a metric~$g_{ab}$, then the Levi-Civita
connection $\widehat\nabla_a$ for the conformally related metric
$\widehat{g}_{ab}=\Omega^2g_{ab}$ is given on $1$-forms by
\begin{gather}\label{change}
\widehat\nabla_a\phi_b=
\nabla_a\phi_b-\Upsilon_a\phi_b-\Upsilon_b\phi_a+\Upsilon^c\phi_cg_{ab},
\end{gather}
where $\Upsilon_a=\nabla_a\log\Omega$ and $\Upsilon^c=g^{cd}\Upsilon_d$. Thus,
associated with a conformal class of metrics we have a restricted
supply of torsion-free connections related by formulae such as
(\ref{contorsion}) where
\begin{gather}\label{thisisGamma}\Gamma_{ab}{}^c=
\Upsilon_a\delta_b{}^c+\Upsilon_b\delta_a{}^c-g_{ab}\Upsilon^c
\end{gather}
for closed $1$-forms $\Upsilon_a$. It is convenient introduce a line bundle $L$
on a general conformal manifold as follows. A choice of metric $g_{ab}$ in the
conformal class trivialises~$L$, allowing us to view a section of $L$ as a
function~$f$. We decree, however, that viewing the same section by means of the
trivialisation due to $\widehat{g}_{ab}=\Omega^2g_{ab}$ gives us the function
$\widehat{f}=\Omega f$. More generally, sections of the bundle $L^w$ transform
by $\widehat{f}=\Omega^wf$ and are called conformal {\em densities} of
conformal {\em weight}~$w$. Since
\[\nabla_a\widehat{f}=\nabla_a(\Omega^wf)=\Omega^w(\nabla_af+w\Upsilon_af)\]
we may view $L^w$ as equipped with a family of connections
related by
\begin{gather}\label{densities}
\widehat\nabla_af=\nabla_af+w\Upsilon_af.
\end{gather}
On an oriented conformal $n$-manifold, the line bundle $\Lambda^n$ is also
trivialised by a choice of metric. Specif\/ically, we may use the inverse metric
$g^{ab}$ to raise indices and then normalise the volume form
$\epsilon_{bc\cdots de}$ so that
$\epsilon^{bc\cdots de}\epsilon_{bc\cdots de}=n!$. Since
$\widehat{g}^{ab}=\Omega^{-2}g^{ab}$ it follows that
$\widehat\epsilon_{bc\cdots de}=\Omega^n\epsilon_{bc\cdots de}$ and hence
that we may identify $\Lambda^n=L^{-n}$. It is easily conf\/irmed that
(\ref{change}) implies
\[\widehat\nabla_a\phi_{bc\cdots de}=
\nabla_a\phi_{bc\cdots de}-n\Upsilon_a\phi_{bc\cdots de}\]
for $n$-forms $\phi_{bc\cdots de}$, which is consistent with (\ref{densities})
for~$w=-n$, as it should be. Combining~(\ref{change}) and~(\ref{densities})
gives us a formula for the change of connection on $\Lambda^1\otimes L^w$,
namely
\[\widehat\nabla_a\phi_b=
\nabla_a\phi_b+(w-1)\Upsilon_a\phi_b-\Upsilon_b\phi_a+\Upsilon^c\phi_cg_{ab}.\]
Let us now consider how the Laplacian $\Delta=\nabla_a\nabla^a$ acting on
densities of weight $w$ is af\/fected by a conformal rescaling of the metric. We
compute:
\[\begin{array}{rcl}
\widehat\nabla_a\widehat\nabla_bf&=&
\nabla_a\widehat\nabla_bf
+(w-1)\Upsilon_a\widehat\nabla_bf
-\Upsilon_b\widehat\nabla_af
+(\Upsilon^c\widehat\nabla_cf)g_{ab}
\end{array}\]
and so
\begin{gather*}
\widehat\Delta f=\widehat{g}^{ab}\widehat\nabla_a\widehat\nabla_bf=
\Omega^{-2}g^{ab}\big(\nabla_a\widehat\nabla_bf
+(w-1)\Upsilon_a\widehat\nabla_bf
-\Upsilon_b\widehat\nabla_af
+(\Upsilon^c\widehat\nabla_cf)g_{ab}\big)\\
\phantom{\widehat\Delta f=\widehat{g}^{ab}\widehat\nabla_a\widehat\nabla_bf}{} =
\Omega^{-2}\big(\nabla^a\widehat\nabla_af
+(n+w-2)\Upsilon^a\widehat\nabla_af\big)\\
\phantom{\widehat\Delta f=\widehat{g}^{ab}\widehat\nabla_a\widehat\nabla_bf}{} =
\Omega^{-2}\big(\nabla^a(\nabla_af+w\Upsilon_af)
+(n+w-2)\Upsilon^a(\nabla_af+w\Upsilon_af)\big).
\end{gather*}
Here, we are using $g^{ab}$ to raise indices on the right hand side. Regrouping
as
\[\widehat\Delta f=\Omega^{-2}\big(\Delta f+(n+2w-2)\Upsilon^a\nabla_af
+w(\nabla^a\Upsilon_a+(n+w-2)\Upsilon^a\Upsilon_a)f\big),\]
we see that when $w=1-n/2$ there are no f\/irst order derivatives in $f$ and
\begin{gather}\label{changeinDelta}
\widehat\Delta f=\Omega^{-2}
\big(\Delta f-(n/2-1)(\nabla^a\Upsilon_a+(n/2-1)\Upsilon^a\Upsilon_a)f\big).
\end{gather}
On the other hand, the Riemann curvature tensor transforms by
\[\widehat R_{abcd}=\Omega^2
(R_{abcd}-\Xi_{ac}g_{bd}+\Xi_{bc}g_{ad}-\Xi_{bd}g_{ac}+\Xi_{ad}g_{bc})
\]
where
\[\Xi_{ab}\equiv\nabla_a\Upsilon_b-\Upsilon_a\Upsilon_b+
{\tfrac{1}{2}}\Upsilon_c\Upsilon^cg_{ab}.\]
It follows immediately that the scalar curvature $R=g^{ac}g^{cd}R_{abcd}$
transforms by
\begin{gather}\label{changeinR}\widehat R=
\Omega^{-2}\big(R-2(n-1)(\nabla^a\Upsilon_a+(n/2-1)\Upsilon^a\Upsilon_a)\big).
\end{gather}
{From} (\ref{changeinDelta}) and (\ref{changeinR}) we conclude that
\[\widehat\Delta f-\tfrac{n-2}{4(n-1)}\widehat Rf=
\Omega^{-2}
\left(\Delta f-\tfrac{n-2}{4(n-1)}Rf\right).\]
If we also absorb the factor of $\Omega^{-2}$ into a change of conformal
weight, then we conclude that $Y\equiv\Delta-\frac{n-2}{4(n-1)}R$ is a
{\em conformally invariant differential operator\/} $Y:L^{1-n/2}\to L^{-1-n/2}$.
This is the {\em conformal Laplacian\/} or {\em Yamabe operator}. The peculiar
multiple of the scalar curvature that must be added to the Laplacian to achieve
conformal invariance is often referred to a {\em curvature correction\/} term.

Now we are in a position to compare with the f\/irst order symmetries of the
Laplacian on~${\mathbb{R}}^n$ found in~Section~\ref{symmetries}. Of course, the
curvature correction terms are absent for the f\/lat metric on~${\mathbb{R}}^n$.
Therefore, there is no dif\/ference between the Laplacian and the Yamabe operator
provided we restrict our attention to conformal rescalings that take the f\/lat
metric to another f\/lat metric. But f\/lat-to-f\/lat rescalings are precisely what
are provided by the conformal Killing f\/ields (\ref{cKillingfield}).
Specif\/ically, the conformal Killing f\/ields are those vector f\/ields whose f\/lows
preserve the conformal structure on~${\mathbb{R}}^n$ and, moreover, the
f\/lat-to-f\/lat conformal factors generated in this way are
general~\cite{bairdwood}. The conformal invariance of the Yamabe operator on a
general conformal manifold therefore implies the invariance of the Laplacian on
${\mathbb{R}}^n$ under Lie derivative ${\mathcal{L}}_V$ for conformal Killing
f\/ields~$V$. Here, as we have just seen, ${\mathcal{L}}_V$ should be interpreted
as acting on the bundles $L^{1-n/2}$ and $L^{-1-n/2}$. At the end of
Section~\ref{symmetries} we found that the symmetry ${\mathcal{D}}_V$ given by
(\ref{DVf}) may be regarded as the Lie derivative on $(\Lambda^n)^{(n-2)/2n}$.
These viewpoints now coincide because
\[\Lambda^n=L^{-n}\implies(\Lambda^n)^{(n-2)/2n}=L^{1-n/2}.\]

The same reasoning applies more generally. The geometric interpretation of
conformal Killing f\/ields combines with the conformal invariance of the Yamabe
operator and we have proved:
\begin{theorem}\label{obvious}
Suppose that $V^a$ is a conformal Killing field on a Riemannian manifold, i.e.\
a~solution of the equation~\eqref{cKillingfield}. Then ${\mathcal{D}}_V$
given by \eqref{DVf} is a symmetry of the Yamabe operator. More precisely,
\[ \left(\Delta -\tfrac{n-2}{4(n-1)}R\right){\mathcal{D}}_V=
\delta_V\left(\Delta -\tfrac{n-2}{4(n-1)}R\right),\]
where
\[ {\mathcal{D}}_V=V^a\nabla_a+\tfrac{n-2}{2n}(\nabla_aV^a)\qquad
\mbox{and}\qquad
\delta_V=V^a\nabla_a+\tfrac{n+2}{2n}(\nabla_aV^a).\]
\end{theorem}

In fact, it follows from proof of \cite[Theorem~1]{laplace} that
these are the only symmetries of the Yamabe operator.

We have already observed that the dif\/ferential pairing ${\mathcal{D}}_Vf$ given
by (\ref{DVf}) is invariant in the very strong sense of being intrinsically
def\/ined on any manifold. Conformal invariance is a weaker statement but one
that is easily verif\/ied directly. {From} (\ref{contorsion}) and
(\ref{thisisGamma}),
\[\widehat\nabla_aV^b=\nabla_aV^b+\Upsilon_aV^b-\Upsilon^bV_a
+\Upsilon_cV^c\delta_a{}^b \ \implies \
\widehat\nabla_aV^a=\nabla_aV^a+n\Upsilon_aV^a\]
whilst for $f$ a density of weight $w$ we have~(\ref{densities}). It follows
immediately that
\begin{gather}\label{firstexample}
\textstyle V^a\nabla_af-\frac{w}{n}(\nabla_aV^a)f\end{gather}
is a conformally invariant dif\/ferential pairing and for $w=1-n/2$ this agrees
with~(\ref{DVf}).

\section{Higher order symmetries}\label{higher}
There are some second order symmetries of the Laplacian on ${\mathbb{R}}^n$
that should be regarded as trivial. If ${\mathcal{D}}=h\Delta$ for some smooth
function~$h$, then $\Delta{\mathcal{D}}=\delta\Delta$ simply by taking
$\delta f=\Delta(hf)$. All symmetries preserve harmonic functions but
${\mathcal{D}}=h\Delta$ does this by dint of annihilating them. Ignoring these
trivial symmetries, the second order symmetries of the Laplacian in
${\mathbb{R}}^3$ were found by Boyer, Kalnins, and Miller~\cite{bkm}. They add
an extra $35$ dimensions to the $11$-dimensional space~(\ref{eleven}). In
\cite{laplace} all higher symmetries were found in all dimensions and, ignoring
the trivial ones, there is an additional f\/inite-dimensional space of second
order ones. More precisely, there are second order symmetries
\begin{gather}\label{DVVf} {\mathcal{D}}_Vf=V^{ab}\nabla_a\nabla_bf
+\tfrac{n}{n+2}(\nabla_aV^{ab})\nabla_bf
+\tfrac{(n-2)n}{4(n+1)(n+2)}(\nabla_a\nabla_bV^{ab})f\end{gather}
for any trace-free symmetric tensor $V^{ab}$ satisfying
\begin{gather}\label{cKillingtensor}
\nabla_aV_{bc}+\nabla_bV_{ca}+\nabla_cV_{ab}=
\tfrac{2}{n+2}
\big(g_{ab}\nabla^dV_{cd}+g_{bc}\nabla^dV_{ad}+g_{ca}\nabla^dV_{bd}\big).
\end{gather}
Such $V^{ab}$ are called {\em conformal Killing tensors} (of valence $2$) and
on ${\mathbb{R}}^n$ form a vector space of dimension
\[\tfrac{(n-1)(n+2)(n+3)(n+4)}{12}.\]

As detailed in~\cite{laplace}, this is the `additional space' of second order
symmetries referred to above. As with f\/irst order symmetries~(\ref{DVf}), the
pairing between valence $2$ symmetric trace-free tensors~$V^{ab}$ and densities
$f$ of weight $1-n/2$ given by the right hand side of (\ref{DVVf}) is
conformally invariant under f\/lat-to-f\/lat conformal rescalings. Unlike the f\/irst
order case, however, the pairing (\ref{DVVf}) is not invariant under general
conformal rescalings of a general Riemannian manifold but, like the Laplacian,
becomes so with the addition of suitable curvature correction terms.
Specif\/ically, for $V^{ab}$ any trace-free symmetric contravariant tensor and
$f$ of weight~$w$
\begin{gather}\label{secondexample}
V^{ab}\nabla_a\nabla_bf-\tfrac{2(w-1)}{n+2}(\nabla_aV^{ab})\nabla_bf
+\tfrac{w(w-1)}{(n+1)(n+2)}(\nabla_a\nabla_bV^{ab})f+
\tfrac{w(n+w)}{(n+1)(n-2)}R_{ab}V^{ab}f,\end{gather}
where $R_{ab}$ is the Ricci tensor, is a
{\em conformally invariant differential pairing\/} on an arbitrary Riemannian
manifold. When $w=1-n/2$, we obtain
\[ {\mathcal{D}}_Vf\equiv V^{ab}\nabla_a\nabla_bf
+\tfrac{n}{n+2}(\nabla_aV^{ab})\nabla_bf
+\tfrac{(n-2)n}{4(n+1)(n+2)}(\nabla_a\nabla_bV^{ab})f
-\tfrac{n+2}{4(n+1)}R_{ab}V^{ab}f\]
as a conformally invariant analogue of the f\/lat symmetry operator~(\ref{DVVf}).
Whether this operator actually provides a symmetry of the Yamabe operator in
general, however, is unclear. By analogy with Theorem~\ref{obvious}, one might
hope that if $V^{ab}$ were a conformal Killing tensor~(\ref{cKillingtensor}),
then it would follow that
\begin{gather}\label{yamabesymmetry}
\left(\Delta -\tfrac{n-2}{4(n-1)}R\right){\mathcal{D}}_V=
\delta_V\left(\Delta -\tfrac{n-2}{4(n-1)}R\right)\end{gather}
for ${\mathcal{D}}_V$ as above and
\[ \delta_Vf\equiv V^{ab}\nabla_a\nabla_bf
+\tfrac{n+4}{n+2}(\nabla_aV^{ab})\nabla_bf
+\tfrac{n+4}{4(n+1)}(\nabla_a\nabla_bV^{ab})f
-\tfrac{n+2}{4(n+1)}R_{ab}V^{ab}f.\]
This is currently unknown. The problem is that there does not seem to be any
useful geometric interpretation of $V^{ab}$ being a conformal Killing tensor in
parallel to $V^a$ being a conformal Killing f\/ield. In principle, the validity
or otherwise of~(\ref{yamabesymmetry}) should boil down to a calculation once
the dif\/ferential consequences of~(\ref{cKillingtensor}) are determined.
Furthermore, the dif\/ferential consequences of~(\ref{cKillingtensor}) can all be
found by prolonging the system, as described in~\cite{bceg} for example.
Nevertheless, it will surely be a dif\/f\/icult calculation.

\section{Further examples of invariant pairings}
The simple conformally invariant pairings we have found so far, namely
(\ref{firstexample}) and (\ref{secondexample}), may be extended by allowing
$V^a$ and $V^{ab}$ to have a general conformal weight as follows. It is
convenient to set
\[
\Phi_{ab}\equiv\frac{1}{n-2}\left(R_{ab}-\tfrac{1}{n}Rg_{ab}\right).\]
\begin{proposition}
Suppose that $V^a$ has conformal weight $v$ and $f$ has conformal weight~$w$.
Then
\begin{gather}\label{firstexampleextended}
(v+n)V^a\nabla_af-w(\nabla_aV^a)f\end{gather}
is conformally invariant. Suppose $V^{ab}$ is trace-free symmetric with
conformal weight $v$ and $f$ has conformal weight~$w$. Then
\begin{gather}
(n+v+2)(n+v+1)V^{ab}\nabla_a\nabla_bf
-2(w-1)(n+v+1)(\nabla_aV^{ab})\nabla_bf\nonumber\\
\qquad{}+w(w-1)(\nabla_a\nabla_bV^{ab})f+w(n+v+w)(n+v+2)\Phi_{ab}V^{ab}f\label{secondexampleextended}
\end{gather}
is conformally invariant.
\end{proposition}

\begin{proof} These are easily verif\/ied using the formulae for conformal
rescaling developed in Section~\ref{conformalgeometry}.\end{proof}

Notice that there are several special cases occurring for particular values of
$v$ and~$w$. The form of (\ref{firstexampleextended}) suggests setting $w=0$
or $v=-n$ and, when we do, we f\/ind
\[(v+n)V^a\nabla_af\quad\mbox{is invariant}\quad\forall\,V^a\implies
f\mapsto\nabla_af\quad\mbox{is invariant}\]
and
\[-w(\nabla_aV^a)f\quad\mbox{is invariant}\quad\forall\,f\implies
V^a\mapsto\nabla_aV^a\quad\mbox{is invariant}.\]
These are familiar invariant linear dif\/ferential operators. The f\/irst one is
the exterior derivative $d:\Lambda^0\to\Lambda^1$ and the second may be
identif\/ied as $d:\Lambda^{n-1}\to\Lambda^n$. More interesting are the
consequences of setting $w=1$ and $v=-n-1$ in (\ref{secondexampleextended}). We
obtain conformally invariant linear dif\/ferential operators
\[
f \mapsto \nabla_a\nabla_bf-\tfrac 1n (\Delta f) g_{ab}+\Phi_{ab}f\qquad\mbox{and}\qquad
V^{ab}\mapsto\nabla_a\nabla_bV^{ab}+\Phi_{ab}V^{ab},\qquad\mbox{respectively}.
\]
The theory of conformally invariant linear dif\/ferential operators is well
understood~\cite{be,b,er,f} and includes these two as simple examples. Other
special values of the weights $v$ and $w$ show up as zeroes of the
coef\/f\/icients in~(\ref{secondexampleextended}). We have already seen that for
$w=0$ there is an invariant dif\/ferential operator $f\mapsto\nabla_af$.
Similarly, for $v=-n-2$ the operator $V^{ab}\mapsto\nabla_aV^{ab}$ is
conformally invariant.

It is also interesting to note that, conversely, we can build the invariant
pairing (\ref{secondexampleextended}) from these various linear invariant
operators. Specif\/ically,
if $V^{ab}$ has weight $v$ and $f$ has weight~$w\not=0$ then, at least where
$f$ does not vanish,
\[f^{-(n+v+1)/w}V^{ab}\mbox{ has weight }-n-1\qquad\mbox{and}\qquad
f^{-(n+v+2)/w}V^{ab}\mbox{ has weight }-n-2.\] Therefore,
\[\nabla_a\nabla_b(f^{-(n+v+1)/w}V^{ab})+\Phi_{ab}f^{-(n+v+1)/w}V^{ab}\qquad
\mbox{and}\qquad\nabla_b(f^{-(n+v+2)/w}V^{ab})\]
are invariant. If the second of these is multiplied by $f^{2/w}$ we obtain a
vector f\/ield of weight $-n$ and conclude that
\[\nabla_a\big(f^{2/w}\nabla_b(f^{-(n+v+2)/w}V^{ab})\big)\]
is also invariant. It follows, therefore, that the combination
\begin{gather*}w(n+v+2)(n+v+w)f^{(n+v+1+w)/w}\big(
\nabla_a\nabla_b(f^{-(n+v+1)/w}V^{ab})+\Phi_{ab}f^{-(n+v+1)/w}V^{ab}
\big)\\
\qquad{}-w(n+v+1)(n+v+1+w) f^{(n+v+w)/w}\big(
\nabla_a\big(f^{2/w}\nabla_b(f^{-(n+v+2)/w}V^{ab})\big)
\big)
\end{gather*}
is invariant. A short calculation reveals that this expression agrees with
(\ref{secondexampleextended}) and, although we supposed $w\not=0$ and also
that $f\not=0$ in carrying out this derivation, the f\/inal conclusion is valid
regardless. These tricks reveal some sort of relationship between invariant
dif\/ferential pairings and invariant linear dif\/ferential operators, at least
when one of the quantities to be paired is a~conformal density.

These tricks are unavailable more generally but there are, nevertheless, many
conformally invariant dif\/ferential pairings. For example, if $V^{ab}$ is
symmetric trace-free with conformal weight~$v$ and~$\phi_a$ has conformal
weight~$w$, then
\[(n+v+2)V^{ab}\nabla_a\phi_b-(w-2)(\nabla_aV^{ab})\phi_b.\]
is invariant. Operators such as this one and (\ref{firstexampleextended}) are
referred to as {\em first order\/}, not only because there are no higher
derivatives involved than f\/irst order but also because there are no cross terms
of the form $\nabla V\bowtie\nabla\phi$ where $\bowtie$ denotes some algebraic
pairing (we shall come back to this point soon). The theory of f\/irst order
invariant dif\/ferential pairings on a conformal manifold, or more generally on a
manifold with AHS structure, is reasonably well understood thanks to recent
work of Kroeske~\cite{k}. Following the approach of Fegan~\cite{f} in the
conformal case and \v{C}ap, Slov\'ak, and Sou\v{c}ek~\cite{cssiii} more
generally, Kroeske classif\/ies the f\/irst order invariant pairings provided
that the weights of the tensors involved are so that invariant linear
dif\/ferential operators are excluded. This corresponds well to the relationship
that seemed to be emerging in our examples. For higher order pairings, however,
Kroeske~\cite{k} f\/inds additional unexpected pairings and is able to present a
satisfactory theory only in the f\/lat projective setting.

\section{Some precise formulations}
When $E$ and $F$ are smooth vector bundles on some smooth manifold, it is
well-known that linear dif\/ferential operators $D:E\to F$ of order $\leq k$ are
in $1$--$1$ correspondence with vector bundle homomorphisms $J^kE\to F$ where
$J^kE$ is the $k^{\mathrm{th}}$ {\em jet bundle\/} of~$E$. Indeed, it is usual
to abuse notation and also write $D:J^kE\to F$, sometimes adopting this
viewpoint as the def\/inition of a~linear dif\/ferential operator. There are
canonical short exact sequences of vector bundles
\begin{gather}\label{jet}
\textstyle 0\to\bigodot^k\!\Lambda^1\otimes E\to J^kE\to J^{k-1}E\to 0,
\end{gather}
where $\bigodot^k\!\Lambda^1$ is the bundle of symmetric covariant tensors of
valence $k$ and the composition
\[\textstyle\bigodot^k\!\Lambda^1\otimes E\to J^kE
\stackrel{D}{\longrightarrow} F\]
is called the {\em symbol} of $D$. For further details, see~\cite{s}
for example.

Dif\/ferential pairings may be formulated similarly. Suppose another vector
bundle $V$ is given and we wish to def\/ine what is a dif\/ferential pairing
$E\times F\to V$. The jet exact sequences (\ref{jet}) may be written as
\[\textstyle J^kE=E\;+\;\Lambda^1\otimes E\;+\;\bigodot^2\!\Lambda^1\otimes E
\;+\;\bigodot^3\!\Lambda^1\otimes E\;+\;\cdots\;+\;
\bigodot^k\!\Lambda^1\otimes E,\]
meaning that $J^kE$ is f\/iltered with successive quotients as shown. It follows
that
\[J^kE\otimes J^kF=E\otimes F\;+\;
\begin{array}c\Lambda^1\otimes E\otimes F\\ \oplus\\
E\otimes\Lambda^1\otimes F\end{array}\;+\;
\begin{array}c\bigodot^2\!\Lambda^1\otimes E\otimes F\\ \oplus\\
\Lambda^1\otimes E\otimes\Lambda^1\otimes F\\ \oplus\\
E\otimes\bigodot^2\!\Lambda^1\otimes F\end{array}\;+\;\cdots
\]
and, following~\cite{k}, we def\/ine the bi-jet bundle $J^k(E,F)$ as the
quotient of $J^kE\otimes J^kF$ corresponding to the f\/irst $k+1$ columns of this
composition series. There are bi-jet exact sequences starting with
\[0\to
\begin{array}c\Lambda^1\otimes E\otimes F\\ \oplus\\
E\otimes\Lambda^1\otimes F\end{array}\to J^1(E,F)\to E\otimes F\to 0\]
and
\[0\to
\begin{array}c\bigodot^2\!\Lambda^1\otimes E\otimes F\\ \oplus\\
\Lambda^1\otimes E\otimes\Lambda^1\otimes F\\ \oplus\\
E\otimes\bigodot^2\!\Lambda^1\otimes F\end{array}\to J^2(E,F)\to J^1(E,F)
\to 0.\]
A dif\/ferential pairing $E\times F\to V$ of order $k$ may now def\/ined as a
homomorphism of vector bundles $J^k(E,F)\to V$. Notice that it is the {\em
total order} of the operator that is constrained to be less than of equal
to~$k$. For f\/irst order pairings, for example, we are excluding terms of the
form $\nabla\psi\bowtie\nabla\phi$ (written in the presence of chosen
connections). For second order pairings we exclude
$\nabla^2\psi\bowtie\nabla\phi$, $\nabla\psi\bowtie\nabla^2\phi$, and
$\nabla^2\psi\bowtie\nabla^2\phi$. The {\em symbol} of a pairing is def\/ined as
the composition
\[\bigoplus_{j=0}^k\textstyle
\bigodot^j\!\Lambda^1\otimes E\otimes\bigodot^{k-j}\!\Lambda^1\otimes F
\enskip\longrightarrow\enskip J^k(E,F)\longrightarrow V.\]
Several other basic notions for linear dif\/ferential operators immediately carry
over to bilinear dif\/ferential pairings. Suppose, for example, that we consider
homogeneous vector bundles on a homogeneous space~$G/P$. If $E$ and
$F$ are induced from $P$-modules ${\mathbb{E}}$ and ${\mathbb{F}}$, then
$J^k(E,F)$ will also be induced from a $P$-module, say
$J^k({\mathbb{E}},{\mathbb{F}})$. If $V$ is induced from~${\mathbb{V}}$, a
$G$-invariant dif\/ferential pairing will correspond to a $P$-module homomorphism
$J^k({\mathbb{E}},{\mathbb{F}})\to{\mathbb{V}}$. In particular, when
$G={\mathrm{SO}}(n+1,1)$ and $G/P$ is the conformal sphere (as detailed in
\cite{srni95} for example) then these $P$-module homomorphisms correspond to
f\/lat-to-f\/lat invariant dif\/ferential pairings. This gives an approach to
classif\/ication as adopted in~\cite{k}. Nevertheless, it remains a dif\/f\/icult
problem.

\section{Other sources of invariant pairings}
Firstly, there are other conformally invariant dif\/ferential operators whose
symmetries give rise to invariant dif\/ferential pairings. The square of the
Laplacian is treated in~\cite{el} and the Dirac operator in~\cite{ess}. There
is also another source of pairings derived from the structure of the
symmetry algebra as follows. Let us consider the symmetries of the Laplacian
as in Sections~\ref{symmetries} and~\ref{higher}. It is clear that the composition
of symmetries is again a symmetry. In particular, if we compose two f\/irst
order symmetries (\ref{DVf}) on ${\mathbb{R}}^n$ then we f\/ind~\cite{laplace}
\[ {\mathcal{D}}_V{\mathcal{D}}_Wf={\mathcal{D}}_{V\circledcirc W}f
+\tfrac{1}{2}{\mathcal{D}}_{[V,W]}f-\tfrac{n-2}{4n(n+1)}{\langle V,W\rangle}f+
\tfrac{1}{n}V^aW_a\Delta f\]
where ${\mathcal{D}}_{V\circledcirc W}$ is given by (\ref{DVVf}) and
\begin{gather*}
(V\circledcirc W)^{ab}=
\tfrac{1}{2}V^aW^b+\tfrac{1}{2}V^bW^a-\tfrac{1}{n}g^{ab}V^cW_c,\\
[V,W]^a=V^b\nabla_bV^a-W^b\nabla_bV^a,\\
\langle V,W\rangle =(n+2)(\nabla_bV^a)(\nabla_aW^b)
 -\tfrac{n+2}{n}(\nabla_aV^a)(\nabla_bW^b)
-\tfrac{n+2}{n}V^a\nabla_a\nabla_bW^b\\
\phantom{\langle V,W\rangle =}{}-\tfrac{n+2}{n}W^a\nabla_a\nabla_bV^b
+\Delta(V_aW^a).
\end{gather*}
The last of these is constant when $V^a$ and $W^a$ are Killing f\/ields on
${\mathbb{R}}^n$ and in this case coincides with their inner product under the
Killing form in ${\mathfrak{so}}(n+1,1)$. Remarkably, there is a curvature
corrected version
\begin{gather*} (n+2)(\nabla_bV^a)(\nabla_aW^b)
-\tfrac{n+2}{n}(\nabla_aV^a)(\nabla_bW^b)
-\tfrac{n+2}{n}V^a\nabla_a\nabla_bW^b\\
\qquad{}-\tfrac{n+2}{n}W^a\nabla_a\nabla_bV^b
+\Delta(V_aW^a)-\tfrac{2(n+2)}{(n-2)}R_{ab}V^aW^b
+\tfrac{2n}{(n-1)(n-2)}R\,V_aW^a\end{gather*} that provides an invariant
dif\/ferential pairing on arbitrary vector f\/ields.

This approach to constructing invariant pairings is sometimes called
{\em quantisation}: the symbol of an operator is specif\/ied, one attempts to
build an invariant operator with this symbol, and then one composes these
operators. Further examples are given by Duval and Ovsienko \cite{do} in the
conformal case and Fox \cite{fox} in the projective case.

Finally, there is a general construction of invariant pairings for very
particular weights called {\em cup products} by Calderbank and
Diemer~\cite{cd}. This construction applies in any {\em parabolic geometry}
but further discussion is beyond the scope of this article.

\subsection*{Acknowledgements}
It is a pleasure to acknowledge many extremely useful conversations with
Vladim\'{\i}r Sou\v{c}ek and Jens Kroeske. Support from the Australian Research
Council is also gratefully acknowledged.

\pdfbookmark[1]{References}{ref}
\LastPageEnding
\end{document}